\documentclass[11pt]{article}
\usepackage{amsmath, amsfonts, amsthm}
\numberwithin{equation}{section}

\begin{document}
\author{Ajai Choudhry}
\Large \title{Expressing three consecutive integers\\
 as sums of three cubes}
\date{}
\maketitle

\begin{abstract} 
This paper is concerned with the problem of expressing three consecutive integers as sums of three cubes. We give several parametric solutions of the problem. We also give somewhat trivial solutions of five or seven consecutive integers that can expressed as sums of three cubes. We conclude the paper with an  open problem regarding four or more consecutive integers expressible as sums of three cubes.
\end{abstract}



\section{Introduction}\label{intro}

Ever since Mordell's observation \cite[p. 505]{Mo2} that he did ``not know anything about the integer solutions of  $X^3+Y^3+Z^3=3$
beyond the existence of the four sets $(1, 1, 1), (4, 4, -5)$ etc.; and it must
be very difficult indeed to find out anything about any other solutions", there has been considerable interest in the representation of integers as a sum of three cubes of integers. 

Remarkable progress has been made in recent years and the following  new representation of the integer 3 as a sum of three cubes was discovered by Booker and Sutherland \cite{BS} in 2020:
\begin{multline}
 569936821221962380720^3+(-569936821113563493509)^3\\
 +(-472715493453327032)^3 =3. \label{Booker}
\end{multline}
Except for integers not expressible as a sum of three cubes because of congruence considerations, representations of all  other  integers $ \leq 100$ are now known  with the representations for $74, 33$ and $42$ having been obtained only in the last seven years (see \cite{Hu}, \cite{Bo} and \cite{BS}).

This paper is concerned with the representation of three  consecutive integers as  sums of three cubes, that is, we need to find integers $n$ such that there is a solution of the  three simultaneous  equations,
\begin{align}
n &= x_1^3+x_2^3+x_3^3, \label{sysn1}\\
n+ 1 & = y_1^3+y_2^3+y_3^3, \label{sysn2}\\
 n+ 2 & = z_1^3+z_2^3+z_3^3, \label{sysn3}
\end{align}
in which  $x_i, y_i, z_i, i=1, 2, 3$, are all integers.

On eliminating $n$ between Eqs. \eqref{sysn1} and \eqref{sysn2}, and again between Eqs. \eqref{sysn2} and \eqref{sysn3}, we get the following two equations, respectively:
\begin{align}
y_1^3+y_2^3+y_3^3 -x_1^3-x_2^3-x_3^3 & =1, \label{eqxy}\\
z_1^3+z_2^3+z_3^3 -y_1^3-y_2^3-y_3^3 & =1. \label{eqyz}
\end{align}

If two of the six integers $y_1, y_2, y_3,  -x_1, -x_2, -x_3$, are 0 and 1, and the remaining integers consist of two pairs such that the sum of the two integers in each pair is 0, we get a trivial solution of Eq. \eqref{eqxy}. If the six integers $z_1, z_2, z_3,- y_1, -y_2, -y_3$, satisfy similar conditions, we get a  trivial solution of  Eq. \eqref{eqyz}.

A solution of the simultaneous Eqs. \eqref{sysn1}, \eqref{sysn2} and  \eqref{sysn3} will be considered trivial if both Eqs. \eqref{eqxy} and \eqref{eqyz} are trivially satisfied. If only one of the two equations \eqref{eqxy} and \eqref{eqyz} is trivially satisfied, the solution will be considered  semi-trivial, and if neither of the  
 two Eqs. \eqref{eqxy} and \eqref{eqyz} is  trivially satisfied,  the solution will be considered a nontrivial solution.

Examples of trivial solutions of the simultaneous Eqs.  \eqref{sysn1}, \eqref{sysn2}  and  \eqref{sysn3} are  given by
\[ n=a^3 +0^3 +0^3,\quad n+1 =a^3+1^3 +0^3,\quad n+2 =a^3+1^3+1^3,\]
and 
\[
n=a^3+b^3+(-1)^3, \quad n+1=a^3+b^3+0^3, \quad n+2=a^3+b^3+1^3,
\]
where $a$ and $b$ are arbitrary integers. 

In the next section we obtain semi-trivial solutions of our problem while in Section \ref{nontrivsols} we  obtain  several nontrivial parametric solutions of the simultaneous Eqs. \eqref{sysn1}, \eqref{sysn2}  and  \eqref{sysn3}. We conclude the paper with some remarks about more than three consecutive integers that are expressible as sums of three cubes and a related open problem. 

\section{Semi-trivial examples of three consecutive integers expressible as sums of three cubes}\label{semitrivsols}
To obtain semi-trivial solutions of Eqs. \eqref{sysn1}, \eqref{sysn2} and \eqref{sysn3}, we start with trivial solutions of Eq. \eqref{eqxy} and then obtain nontrivial solutions of the resulting equation obtained from  Eq. \eqref{eqyz}. If we start  with trivial solutions of Eq. \eqref{eqyz} and then solve Eq. \eqref{eqxy}, we get similar results.

As is readily seen, trivial solutions of Eq. \eqref{eqxy} may be written, without loss of generality, in one of the following six ways:
\begin{align}
 (x_1, x_2, x_3, y_1, y_2, y_3)  & = (0, u, v, 1, u, v), \label{trivsol1} \\
(x_1, x_2, x_3, y_1, y_2, y_3)  & = (0, u, -u, 1, v, -v), \label{trivsol2} \\
(x_1, x_2, x_3, y_1, y_2, y_3)  & = ( 0, -1,u, u, v, -v), \label{trivsol3} \\
(x_1, x_2, x_3, y_1, y_2, y_3)  & = (-1, u, v, 0, u, v), \label{trivsol4} \\
(x_1, x_2, x_3, y_1, y_2, y_3)  & = (-1, u, -u, 0, v, -v), \label{trivsol5} \\
(x_1, x_2, x_3, y_1, y_2, y_3)  & = (u, v, -v,  0, 1,u), \label{trivsol6}
\end{align}
where $u$ and $v$ are arbitrary integer parameters.

It follows from the above trivial solutions  of Eq. \eqref{eqxy} that at least two of the nine integers $x_i, y_i, z_i, i=1, 2, 3$, in any  semi-trivial solution of the simultaneous diophantine Eqs. \eqref{sysn1}, \eqref{sysn2}  and  \eqref{sysn3} will necessarily be 0 and 1.

On substituting the values of $x_i, y_i$ given by the six  trivial solutions \eqref{trivsol1}--\eqref{trivsol6},  in succession, in Eq. \eqref{eqyz}, we get the following six equations, respectively:
\begin{align}
-u^3 - v^3 + z_1^3 + z_2^3 + z_3^3 = 2, \label{trivsol1a} \\
z_1^3 + z_2^3 + z_3^3 = 2, \label{trivsol2a} \\
-u^3 + z_1^3 + z_2^3 + z_3^3 = 1, \label{trivsol3a} \\
-u^3 - v^3 + z_1^3 + z_2^3 + z_3^3 = 1, \label{trivsol4a} \\
z_1^3 + z_2^3 + z_3^3 = 1, \label{trivsol5a}\\
-u^3 + z_1^3 + z_2^3 + z_3^3 = 2. \label{trivsol6a}    
\end{align}
 
Nontrivial identities expressing the integers 1 and 2 as a sum of three cubes of polynomials have been given by Mahler \cite{Ma} and by Werebrusow (as stated by Mordell \cite{Mo1}), respectively.  Further, Choudhry \cite{Ch} has given identities expressing the integers 1 and 2 as a sum of four  cubes of polynomials and has also described a method of obtaining identities expressing any arbitrary integer as a sum of five cubes of polynomials, and thus the integers 1 and 2 can also be expressed in this manner. 

Using the  identities expressing 1 or 2 as a sum of cubes of three, four or five polynomials, we immediately get nontrivial parametric solutions of the six equations \eqref{trivsol1a}--\eqref{trivsol6a}. We can thus obtain several parametric  solutions of the simultaneous diophantine Eqs. \eqref{eqxy} and   \eqref{eqyz} in which Eq. \eqref{eqxy} is trivially satisfied,   and we thus obtain semi-trivial parametric solutions  of the three simultaneous  Eqs. \eqref{sysn1}, \eqref{sysn2} and \eqref{sysn3}.

As an example, using the identity
\begin{equation}
(t^2)^3+(t^2)^3+(-t^2 + t + 1)^3+(-t^2 - t + 1)^3=2, \label{identint2first}
\end{equation}
given by Choudhry \cite[p. 4]{Ch},
we get the following solution of \eqref{trivsol6a}:
\begin{equation}
 u=t^2 + t - 1, z_1=t^2, z_2=t^2, z_3=-t^2 + t + 1. \label{soltrivsol6a}
\end{equation}

The solution \eqref{soltrivsol6a} yields  the following solution of the simultaneous Eqs. \eqref{eqxy} and  \eqref{eqyz}:
\begin{equation*}
\begin{aligned}
x_1 &= t^2 + t - 1,  & x_2 &= v,  & x_3 &= -v,  & y_1 &= 0,\;\;  y_2 = 1, \\
y_3 &= t^2 + t - 1,  & z_1 &= t^2,  & z_2 &= t^2,  & z_3 &= -t^2 + t + 1,
\end{aligned}
\end{equation*}
where $t$ and $v$ are arbitrary integer parameters, and hence we obtain the integer $n=(t^2 + t - 1)^3$ such that the three consecutive integers $n, n+1$, and $n+2$ are expressible as $x_1^3+x_2^3+x_3^3, y_1^3+y_2^3+y_3^3, z_1^3+z_2^3+z_3^3$, respectively. Taking $(t,v)=(2,1)$, yields the numerical example,
\[
125 = 5^3+1^3+(-1)^3, 126= 0^3+1^3+5^3, 127 =4^3+4^3+(-1)^3.
\]

\section{Nontrivial examples of three consecutive integers expressible as sums of three cubes}\label{nontrivsols}
We will describe, in the next two subsections, two different methods of obtaining nontrivial solutions of the simultaneous diophantine Eqs. \eqref{sysn1}, \eqref{sysn2} and \eqref{sysn3}, or the equivalent pair of simultaneous Eqs. \eqref{eqxy} and \eqref{eqyz}. The first method generates several multi-parameter  solutions in polynomials of high degree and some of these solutions yield numerical solutions in positive integers. The second method generates  several solutions in terms of linear and quadratic polynomials but all such solutions that we could obtain yield numerical examples that necessarily include negative integers
\subsection{First method}\label{nontrivsol1}
We will first solve the simultaneous  Eqs.  \eqref{eqxy} and \eqref{eqyz} by rewriting  Eq.\eqref{eqxy} as
\begin{equation}
y_1^3+y_2^3+y_3^3 -x_3^3  =x_1^3+x_2^3 +1.  \label{eqxyvar}
\end{equation}
and solving Eq.\eqref{eqxyvar} together with the following equation obtained by eliminating $y_1, y_2, y_3$, from Eqs. \eqref{eqxy} and \eqref{eqyz}:
\begin{equation}
 z_1^3 + z_2^3 + z_3^3 - x_1^3 - x_2^3 - x_3^3 = 2. \label{eqzx}
\end{equation}

We now impose the auxiliary condition $z_3=x_3$ when Eq. \eqref{eqzx} reduces to
\begin{equation}
 z_1^3 + z_2^3  - x_1^3 - x_2^3  = 2. \label{eqzxred}
\end{equation}
While the auxiliary condition $z_3=x_3$ facilitates the solution of our problem, it also ensures that the solution obtained by this method will never consist of distinct integers.

To solve Eq. \eqref{eqxyvar},  we may  use the values of $x_1, x_2$ given by any parametric solution of Eq. \eqref{eqzxred}. One such solution, that immediately follows from the identity \eqref{identint2first}, may be written, in terms of an arbitrary parameter $t$, as follows:
\begin{equation}
x_1=-t^2, \quad x_2=t^2-t-1, \quad z_1=t^2, \quad z_2=-t^2 - t + 1, \label{eqzxsol1}.
\end{equation}
A second solution of Eq. \eqref{eqzxred} is as follows: 
\begin{equation}
\begin{aligned}
x_1 & = 6gt^2(g^3 + h^3), & x_2 & = 6t^3(g^3 + h^3)^2 - 1, \\
 z_1 & = 6t^3(g^3 + h^3)^2 + 1, & z_2 & = -6ht^2(g^3 + h^3)
\end{aligned}
\label{eqzxsol2}
\end{equation}
where $g, h$ and $t$ are arbitrary parameters. The solution \eqref{eqzxsol2} follows from an identity, given by Choudhry \cite[p.\ 4]{Ch}, expressing 2 as a sum of four cubes of polynomials in three variables. It may also be verified by direct computation.

With the chosen values of $x_1, x_2$, we will solve Eq. \eqref{eqxyvar} by assigning values to the four variables $y_1, y_2, y_3, x_3$ in terms of certain parameters such that the left-hand side of Eq. \eqref{eqxyvar} reduces to a polynomial having at least one of the parameters in degree 1, and then Eq. \eqref{eqxyvar} can be readily solved. We describe two ways in which this can be done. The first way is by writing
\begin{equation}
 \begin{aligned}
x_3 & = m + p - q, \quad & y_1 & = m + p + q, \\
  y_2 & = -m + p - q, \quad & y_3 & = m - p - q, 
\end{aligned}
\label{subsxysimple}
\end{equation}  
where $m, p$ and $q$ are arbitrary parameters, and now  the left-hand side of Eq. \eqref{eqxyvar} reduces to $24mpq$. 

The second way is to write
\begin{equation}
x_3=um+p,  y_1=v_1m+p,  y_2=v_2m+q,  y_3=v_3m-q, \label{subsxygen}
\end{equation} 
where $u$ and $v_i, i=1, 2, 3$, are  chosen so as to satisfy the condition,
\begin{equation}
u^3 = v_1^3 + v_2^3 + v_3^3, \label{conduv}
\end{equation}
while $m, p$ and $q$ are arbitrary parameters. Now the left-hand side of Eq. \eqref{eqxyvar} reduces to 
\begin{equation*}
3((v_1^2-u^2)p + (v_2^2 - v_3^2)q)m^2 - 3((u-v_1)p^2- (v_2 + v_3)q^2)m,
\end{equation*}
and, on   choosing $p=(v_2^2 - v_3^2), q=u^2-v_1^2$, it  further reduces to
\begin{multline*}
\quad \quad \quad \quad 3m(u - v_1)(v_2 + v_3)(u^3 + u^2v_1 - uv_1^2 - v_1^3 \\
- v_2^3+ v_2^2v_3 + v_2v_3^2 - v_3^3), \quad \quad \quad \quad 
\end{multline*}
where $m$ occurs only in degree 1. 

We can use  any set of values of $x_1, x_2$ satisfying Eq. \eqref{eqzxred} with either of the two ways given above of reducing the left-hand side of Eq. \eqref{eqxyvar} to a linear polynomial in the independent parameter $m$ and thus obtain several parametric solutions of the simultaneous diophantine Eqs. \eqref{eqxy} and \eqref{eqyz}.

As an example, taking the values of $x_1, x_2$ given by \eqref{eqzxsol1}, and the values of $x_3, y_1,. y_2, y_3$, given by \eqref{subsxysimple}, Eq. \eqref{eqxyvar} reduces to 
\begin{equation}
24mpq =  -t(3t^4 - 5t^2 + 3), \label{eqxyvarred}
\end{equation}
and and we can readily find several parametric solutions by a suitable choice of $t$. For instance, we obtain a three-parameter solution by taking 
\begin{equation}
t=-24pqr, \quad  m=3r(331776p^4q^4r^4 - 960p^2q^2r^2 + 1).
\end{equation}
This yields the following nontrivial solution of the simultaneous Eqs. \eqref{eqxy} and \eqref{eqyz}:
\begin{equation}
\begin{aligned}
x_1 & = -576p^2q^2r^2, \\
x_2 & = 576p^2q^2r^2 + 24pqr - 1, \\
x_3 & = 995328p^4q^4r^5 - 2880p^2q^2r^3 + p - q + 3r, \\
y_1 & = 995328p^4q^4r^5 - 2880p^2q^2r^3 + p + q + 3r, \\
y_2 & = -995328p^4q^4r^5 + 2880p^2q^2r^3 + p - q - 3r, \\
y_3 & = 995328p^4q^4r^5 - 2880p^2q^2r^3 - p - q + 3r, \\
z_1 & = 576p^2q^2r^2, \\
z_2 & = -576p^2q^2r^2 + 24pqr + 1, \\
z_3 & = 995328p^4q^4r^5 - 2880p^2q^2r^3 + p - q + 3r,
\end{aligned}
\label{soleqxyeqyznt1}
\end{equation}
where $p, q$ and $r$ are arbitrary integer parameters.

We can obtain another solution of the simultaneous Eqs. \eqref{eqxy} and \eqref{eqyz} by following a similar procedure  using the values of $x_1, x_2$ given by the second  solution \eqref{eqzxsol2} of Eq. \eqref{eqzxred}. As it is cumbersome to write this solution, we do not give it explicitly.

We will  now obtain solutions of the simultaneous Eqs. \eqref{eqxy} and \eqref{eqyz} by following  the second way mentioned above using the values of $x_3, y_1, y_2, y_3$, given by \eqref{subsxygen}. We can use the values of $x_1, x_2$ given by either of the two solutions \eqref{eqzxsol1} and \eqref{eqzxsol2} of Eq. \eqref{eqzxred}, but we restrict ourselves only  to  the second  solution \eqref{eqzxsol2} of Eq. \eqref{eqzxred}.

We need values of $u, v_1, v_2, v_3$, satisfying the condition \eqref{conduv}. Several parametric solutions of Eq. \eqref{conduv} are already known (\cite[p.\ 257--260]{HW}, \cite[p.\ 290--291]{Hua}, and   starting with such solutions of Eq. \eqref{conduv}, we can obtain multi-parameter solutions of the simultaneous Eqs. \eqref{eqxy} and \eqref{eqyz}.

As the more general multi-parameter solutions obtained by this method are too cumbersome to write, we give  a simpler example taking $(u, v_1, v_2, v_3)=(9, 6, 1, 8)$  when the condition \eqref{conduv} is satisfied. We use the relations \eqref{subsxygen} and take $(p, q) = (7, -5)$, and further, we take  $t=13r$ in the solution \eqref{eqzxsol2}, when Eq. \eqref{eqxyvar} simply reduces to
\begin{multline}
m=169r^3(g + h)^2(g^2 - gh + h^2)^2(57921708g^{12}r^6\\
 + 231686832g^9h^3r^6+ 347530248g^6h^6r^6 + 231686832g^3h^9r^6\\
 + 57921708h^{12}r^6 + 13182g^6r^3 - 13182h^6r^3 + 1). \label{valmsec}
\end{multline}

Thus a solution of the  simultaneous Eqs. \eqref{eqxy} and \eqref{eqyz} is given by \eqref{subsxygen} where $(u, v_1, v_2, v_3, p, q) =(9, 6, 1, 8,$ $7, -5)$ and $m$ is given by \eqref{valmsec} in terms of three arbitrary integer parameters  $g, h$ and $r$. As a  special case we take $(g, h)=(2, -1)$, when we get the following solution of the simultaneous Eqs. \eqref{eqxy} and \eqref{eqyz}:
\begin{equation}
\begin{aligned}
x_1 & =  14196r^2, \\
x_2  & =  645918r^3 - 1, \\
x_3 & =  10364749588252332r^9 + 61893800514r^6 \\ & \quad \quad  + 74529r^3 + 7,\\
 y_1 & =  6909833058834888r^9 + 41262533676r^6 \\ & \quad \quad + 49686r^3 + 7,\\
 y_2 & =  1151638843139148r^9 + 6877088946r^6 \\ & \quad \quad + 8281r^3 - 5, \\
y_3 & =  9213110745113184r^9 + 55016711568r^6 \\ & \quad \quad + 66248r^3 + 5, \\
z_1 & =  645918r^3 + 1, \\
z_2 & =  7098r^2, \\
z_3 & =  10364749588252332r^9 + 61893800514r^6 \\ & \quad \quad + 74529r^3 + 7,
\end{aligned}
\label{soleqxyeqyznt2}
\end{equation}
where $r$ is an arbitrary integer parameter.

We now have two parametric solutions \eqref{soleqxyeqyznt1} and \eqref{soleqxyeqyznt2} of the simultaneous diophantine  Eqs. \eqref{eqxy} and \eqref{eqyz}.  With the values of $x_i, y_i, z_i, i=1, 2, 3$, given by  \eqref{soleqxyeqyznt1} and by \eqref{soleqxyeqyznt2}, we may take $n=x_1^3+x_2^3+x_3^3$, and obtain two nontrivial examples of three consecutive integers $n, n+1, n+2$ that may be expressed as $x_1^3+x_2^3+x_3^3$, $y_1^3+y_2^3+y_3^3$,  and $z_1^3+z_2^3+z_3^3$, respectively. 

As a numerical example, taking $(p, q, r)=(2, 1, 1)$ in the solution \eqref{soleqxyeqyznt1}, we get $n=4030102758035382018255$, and the three consecutive integers beginning with $n$ may be written as sums of three cubes as follows:
\begin{align*}
n & = (-2304)^3+2351^3+ 15913732^3, \\
n+1 & = 15913734^3+ (-15913730)^3+15913728^3, \\
n +2 &  = 2304^3+(-2255)^3+ 15913732^3.
\end{align*}

The second solution \eqref{soleqxyeqyznt2}  yields, on taking $r \geq 1$,  infinitely many nontrivial solutions of our problem in positive integers. For instance,  taking $r=1$ in the solution \eqref{soleqxyeqyznt2}, we get 
\[
n= 1113484618981001668543451628004732068607126098717,
\]
and the representations of the three consecutive integers as sums of three cubes of positive integers are  as follows:
\begin{equation*}
 \begin{aligned}
n & = 14196^3 + 645917^3+ 10364811482127382^3,\\
n+1 & = 6909874321418257^3 + 1151645720236370^3 \\
 & \quad \quad + 9213165761891005^3,\\
n+2 & =7098^3+  645919^3+ 10364811482127382^3.
\end{aligned}
\end{equation*}

\subsection{Second method}
We will now obtain  nontrivial solutions of the simultaneous Eqs. \eqref{eqxy} and \eqref{eqyz} by first obtaining  parametric solutions of  the corresponding simultaneous homogeneous equations namely,
\begin{align}
y_1^3+y_2^3+y_3^3 -x_1^3-x_2^3-x_3^3 & =t^3, \label{eqxyhom}\\
z_1^3+z_2^3+z_3^3 -y_1^3-y_2^3-y_3^3 & =t^3, \label{eqyzhom}
\end{align}
and then choosing the parameters such that we get $t=1$.

On writing
\begin{equation}
\begin{aligned}
x_1& =  a_1 u+ b_1 v, \quad & x_2& = - a_1 u- b_4 v, \quad & x_3& = - a_4 u- b_1 v, \\
y_1& =  a_2 u+ b_2 v,  \quad &y_2& = - a_2 u,       \quad &y_3& = - b_2 v,  \\
z_1& =  a_3 u+ b_3 v,   \quad &z_2& = - a_3 u+ b_4 v,  \quad & z_3& =  a_4 u- b_3 v,\\
t& =  a_4 u+ b_4 v,
\end{aligned}
\label{subsxynt2}
\end{equation}
where $u, v, a_i, b_i, i=1, \ldots, 4$, are all arbitrary parameters, Eqs. \eqref{eqxyhom} and \eqref{eqyzhom} reduce, after transposing all terms to one side and removing the factor $3uv$ in both cases,  to the following two equations, respectively:
\begin{multline}
((a_1^2 - a_4^2)b_1 - a_2^2b_2 - (a_1^2 - a_4^2)b_4)u \\
+ ((a_1 - a_4)b_1^2 - a_2b_2^2 - (a_1 - a_4)b_4^2)v, \label{eqxyhom1}
\end{multline}
and 
\begin{multline}
(a_2^2b_2 - (a_3^2 - a_4^2)b_3 - (a_3^2 - a_4^2)b_4)u\\
 + (a_2b_2^2 - (a_3 + a_4)b_3^2 + (a_3 + a_4)b_4^2)v. \label{eqyzhom1}
\end{multline}

We now equate to 0 the coefficients of $u$ and $v$ in Eq. \eqref{eqxyhom1} and the coefficient of $u$ in Eq. \eqref{eqyzhom1} and solve for $b_i, i=1, \ldots, 4$ excluding   solutions in which $b_2=0$ since that would make $y_3=0$. We thus get the following solution for $b_i, i=1, \ldots, 4$:
\begin{equation}
\begin{aligned}
b_1& = k(a_3^2 - a_4^2)(a_1^3 + a_1^2a_4 - a_1a_4^2 + a_2^3 - a_4^3),\\
b_2& = 2ka_2(a_1^2-a_4^2)(a_3^2-a_4^2),\\
 b_3& = -k((a_3^2-a_4^2)a_1^3-(2a_2^3 - a_3^2a_4 + a_4^3)a_1^2 \\
& \quad \quad - a_4^2(a_3^2 - a_4^2)a_1 - a_2^3a_3^2 + 3a_2^3a_4^2 - a_3^2a_4^3 + a_4^5),\\
 b_4& = k(a_3^2-a_4^2)(a_1^3 + a_1^2a_4 - a_1a_4^2 - a_2^3 - a_4^3),
\end{aligned}
\label{valb1234}
\end{equation}
where $k$ is an arbitrary parameter.

With the values of $b_i, i=1, \ldots, 4$, given by \eqref{valb1234}, Eq. \eqref{eqxyhom1} is identically satisfied for all $u$ and $v$ while Eq. \eqref{eqyzhom1} reduces to the following equation:
\begin{multline}
(a_ 3^2 - a_ 4^2) a_ 1^3 - (a_ 2^3 - a_ 3^3 + a_ 3 a_ 4^2) a_ 1^2 \\
- (a_ 3^2 - a_ 4^2) a_ 4^2a_ 1 - (a_ 3^2 - 2 a_ 4^2) a_ 2^3 - a_ 3^3 a_ 4^2 + a_ 3 a_ 4^4. \label{conda}
\end{multline}

Eq. \eqref{conda} has a parametric solution given by
\begin{equation}
\begin{aligned}
a_1 & =  -p^2 - p  q + q^2, & a_2 & =  2  p  q, \\
a_3 & =  p^2 - p  q - q^2,   & a_4 & =  p^2 + p  q + q^2,
\end{aligned} \label{parmsola}
\end{equation}
where $p$ and $q$ are arbitrary parameters. It is also possible to find infinitely many integer solutions of Eq. \eqref{conda} that are not given by the parametric solution \eqref{parmsola}. While these solutions of Eq. \eqref{conda} will yield solutions of the simultaneous diophantine Eqs.  \eqref{eqxyhom} and \eqref{eqyzhom} in terms of two linear parameters $u$ and $v$, in order to obtain solutions of the simultaneous diophantine Eqs.  \eqref{eqxy} and \eqref{eqyz},  we need to find such solutions in which we can find $u, v$ so as to satisfy the condition $t=a_4 u+ b_4 v=1$. For this purpose, it suffices to find integer solutions of Eq. \eqref{conda} such that the resulting  value of $b_4$, obtained from \eqref{valb1234}, is coprime with $a_4$. 

Accordingly, we performed trials over the range $|a_1| +|a_2| +|a_3| + |a_4| \leq 100$, and found several solutions of Eq. \eqref{conda} such that $\gcd (a_4, b_4)=1$. However, only two of these solutions yielded independent solutions of the simultaneous diophantine Eqs.  \eqref{eqxy} and \eqref{eqyz}. 

One  of the aforementioned two solutions of  Eq. \eqref{conda} namely,  $(a_1, a_2, a_3, a_4) = (5, 2, -2, 1)$, yields on using the relations \eqref{valb1234} and taking $k=1/24$, $(b_1, b_2, b_3, b_4)=( 19,  12,$  $-1,  17)$, and now, on using the relations \eqref{subsxynt2}, we get   the following solution of the simultaneous diophantine Eqs.  \eqref{eqxyhom} and \eqref{eqyzhom}:
\begin{equation}
\begin{aligned}
 x_1  & =  5 u + 19 v,   & x_2  & =  -5 u - 17 v,   & x_3  & =  -u - 19 v,   \\
 y_1  & =  2 u + 12 v,   & y_2  & =  -2 u, & y_3  & =  -12 v, \\
z_1  & =  -2 u - v,   & z_2  & =  2 u + 17 v,   & z_3  & =  u + v, \\  
 t  & =  u + 17 v, 
\end{aligned}
\label{sol1hom}
\end{equation}
where $u$ and $v$ are arbitrary parameters. Now, on writing $u=1-17v$, we get $t=1$, and hence  we get the following solution of the simultaneous diophantine Eqs.  \eqref{eqxy} and \eqref{eqyz}:

\begin{equation}
\begin{aligned}
x_1 & =  - 66v +5,  &  x_2 & =   68v -5,  &  x_3 & =   - 2v-1, \\
y_1 & = - 22v+2,  &  y_2 & =   34v-2,  &  y_3 & =  -12v, \\
z_1 & =   33v-2,  &  z_2 & = - 17v+2,  &  z_3 & =   - 16v+1,
\end{aligned}
\label{sol1nt2}
\end{equation}
where $v$ is an arbitrary integer parameter.

With the values of $x_i, y_i, z_i, i=1, 2, 3$, given by  \eqref{sol1nt2}, we get 
\[
n=x_1^3+x_2^3+x_3^3 = 26928v^3 - 4032v^2 + 144v - 1,
\]
and the  three consecutive integers $n, n+1, n+2$  may be expressed as $x_1^3+x_2^3+x_3^3$, $y_1^3+y_2^3+y_3^3$,  and $z_1^3+z_2^3+z_3^3$, respectively.

On taking $v=2$ in the solution \eqref{sol1nt2}, we get the following numerical example of three consecutive integers expresible as sums of three cubes:
\[
\begin{aligned}
199583  & =  131^3+(-5)^3+(-127)^3, \\
199584  & =  66^3+(-24)^3+(-42)^3, \\
199585  & =  64^3+(-31)^3+(-32)^3,
\end{aligned}
\]

A second solution of  Eq. \eqref{conda} namely,  $(a_1, a_2, a_3, a_4) = (-4, 6,$ $19, 1)$, yields, on following a similar procedure as above,   the following solution of the simultaneous diophantine Eqs.  \eqref{eqxy} and \eqref{eqyz}:
\begin{equation}
\begin{aligned}
x_1   & =    -97v - 4,  &    x_2   & =    145v + 4,  &    x_3   & =    -48v - 1, \\
y_1   & =    194v + 6,  &    y_2   & =    -174v - 6,  &    y_3   & =    -20v, \\
z_1   & =    582v + 19,  &    z_2   & =    -580v - 19,  &    z_3   & =    -2v + 1,
\end{aligned}
\label{sol2nt2}
\end{equation}
where, as before,  $v$ is an arbitrary integer parameter. This gives 
\[
n=x_1^3+x_2^3+x_3^3 =2025360  v^3 + 132480  v^2 + 2160  v - 1,
\]
and, as before,  the  three consecutive integers $n, n+1, n+2$  may be expressed as $x_1^3+x_2^3+x_3^3$, $y_1^3+y_2^3+y_3^3$,  and $z_1^3+z_2^3+z_3^3$, respectively.

While the two linear parametric solutions \eqref{sol1nt2} and \eqref{sol2nt2} of the simultaneous diophantine Eqs.  \eqref{eqxy} and \eqref{eqyz} give numerical solutions of these equations in distinct integers, it is easily seen that some of these integers will necessarily be negative. We, accordingly, explored the existence of parametric solutions of degree 2 and obtained one such solution which is as follows:
\begin{equation}
\begin{aligned}
 x_1 & =  -147m^2 - 42m - 1,  & x_2 & =  294m^2 + 77m + 3, \\
 x_3 & =  -147m^2 - 35m - 3, & y_1 & = -147m^2 - 56m - 4, \\
 y_2 & = 294m^2 + 77m + 4, & y_3 & =  -147m^2 - 21m, \\
 z_1 & =  -147m^2 - 14m + 1,&  z_2 & =294m^2 + 77m + 5,\\
z_3 & =  -147m^2 - 63m - 5,  &
\end{aligned}
\label{sol3nt2}
\end{equation}
where $m$ is an arbitrary integer parameter.

The procedure for obtaining the solution \eqref{sol3nt2} is  similar to that used for obtaining the two linear solutions but we omit the tedious details. The solution \eqref{sol3nt2}   may be readily verified by direct computation.

As in the case of earlier solutions, it  follows from \eqref{sysn1} that  
\begin{align*}
n & = 19059138m^6 + 14975037m^5 + 4429845m^4 \\
& \quad \quad+ 617400m^3+ 40572m^2 + 1008m - 1,
\end{align*}
and the relations
\begin{equation*}
n= x_1^3+x_2^3+x_3^3, n+1=y_1^3+y_2^3+y_3^3, n+2=z_1^3+z_2^3+z_3^3,
\end{equation*}
show that the three consecutive integers $n, n+1$ and $n+2$ are all expressible as sums of cubes of integers

As a numerical example, taking $m=1$ in the solution \eqref{sol3nt2}, we get $n=39122999$, and the three consecutive integers beginning with $n$ may be written as  sums of three cubes as follows:
\begin{align*}
39122999 & = 374^3+(-185)^3+(-190)^3, \\
39123000 & = 375^3+(-168)^3+(-207)^3, \\
39123001 & = 376^3+(-160)^3+(-215)^3.
\end{align*}

While the solution \eqref{sol3nt2}  yields numerical solutions of  Eqs. \eqref{eqxy} and \eqref{eqyz} in distinct integers, here also   some of the integers are necessarily negative.

\section{Four or more consecutive integers expressible as sums of three cubes and a related open problem} 
If $a$ is any arbitrary integer, the representations  $a^3, a^3 \pm 1^3 \pm 1^3$ furnish a trivial example of five consecutive integers $a^3-2, \ldots, a^3+2$, that can be expressed as a sum of three cubes.
Less trivial examples are obtained from integer solutions of the diophantine Eq. \eqref{eqzxred}
and taking $n=x_1^3+x_2^3-1$ when the five integers $n, n+1, \ldots, n+4$ are all expressible as sums of three cubes since we have,
\[
\begin{aligned}
n & =x_1^3+x_2^3+(-1)^3, & n+1 &= x_1^3+x_2^3+0^3, \\
n+2 & =  x_1^3+x_2^3+1^3, & n+3 & = z_1^3+z_2^3 +0^3, \\
n+4 & = z_1^3+z_2^3 +1^3. & &
\end{aligned}
\]
Thus the parametric  solutions \eqref{eqzxsol1} and \eqref{eqzxsol2} of Eq. \eqref{eqzxred} yield infinitely many  examples of five consecutive integers expressible as  sums of three cubes.

Further, if the integers $x_i, i=1, 2, 3$, satisfy the relation $x_1^3+x_2^3-x_3^3=3$ and we take $n=x_3^3-2$, then the 7 consecutive integers $n +j, j=0, \ldots, 6$,  can all be represented  as a sum of three cubes since we have
\[
\begin{aligned}
n & =x_3^3 +(-1)^3 +(-1)^3, & n+1 &= x_3^3 +(-1)^3 +0^3, \\
 n+2 &=x_3^3 +0^3 +0^3, &n+3 & = x_1^3+1^3+0^3, \\
 n+4 & =  x_1^3+1^3 +1^3, & n+5 & =x_1^3+x_2^3 +0^3, \\
n+6 & = x_1^3+x_2^3 +1^3. & & 
\end{aligned}
\]

Since $4^3+4^3-5^3=3$, taking $n=5^3-2=123$ gives 7 consecutive integers, commencing with 123,  that  can all be represented  as a sum of three cubes. Similarly, in view of the identity \eqref{Booker}, the 7 consecutive integers starting with $472715493453327032^3-2$ can all be expressed as sums of three cubes.

None of the above  examples of 5 or 7 consecutive integers is really  nontrivial since the representations of several of the integers as a sum of three cubes include the cubes of 0, 1 and $-1$. It is an open problem  of considerable interest to find infinitely many nontrivial examples of four or more consecutive integers that can all be expressed as a sum of three cubes without using the cubes of the integers $0, -1$ and $ 1$.

\noindent Postal Address: Ajai Choudhry, 
\newline \hspace*{1.30 in} 13/4 A Clay Square,
\newline \hspace*{1.30 in} Lucknow - 226001, INDIA.
\newline \noindent  E-mail: ajaic203@yahoo.com

\end{document}